\documentclass[a4paper]{article}
\usepackage{graphicx}
\usepackage{amsmath,amssymb,amsthm,bbm}
\usepackage{todonotes}
\setuptodonotes{fancyline,color=green!50}
\usepackage{url}
\usepackage{tikz}
\usetikzlibrary{calc,shapes,positioning}

\newtheorem{example}{Example}

\newcommand{\RR}{\mathbbm{R}}

\newcommand{\findet}{\Delta_{\mathrm{fin}}}
\newcommand{\infdet}{\Delta_{\mathrm{inf}}}

\begin{document}

\title{Structural Identifiability and Discrete Symmetries}
\author{Xabier Rey Barreiro\thanks{X. Rey Barreiro and A.F. Villaverde are with CITMAga, Santiago de Compostela, Galicia, 15782 Spain, and with the  Universidade de Vigo, Vigo, Galicia, 36310 Spain (e-mails: \{xabier.rey, afvillaverde\}@uvigo.gal).}\ \thanks{Xabier Rey Barreiro and Nick Baberuxki should be considered joint first authors.} \and Nick Baberuxki\thanks{N. Baberuxki, M. Abebaw Mebratie and W.M. Seiler are with Universität Kassel (e-mails:
\{nbaberuxki, meskerem.mebratie, seiler\}@mathematik.uni-kassel.de).}\ \footnotemark[2] \and Meskerem Abebaw Mebratie\footnotemark[3] \and Alejandro F. Villaverde\footnotemark[1]\ \thanks{Corresponding authors.} \and Werner M. Seiler\footnotemark[3]\ \footnotemark[4]}

\date{\today}
\maketitle

\begin{abstract}
  We discuss the use of symmetries for analyzing the structural
  identifiability and observability of control systems.  Special emphasis
  is put on the role of discrete symmetries, in contrast to the more
  commonly studied continuous or Lie symmetries.  We argue that discrete
  symmetries are the origin of parameters which are structurally locally
  identifiable, but not globally.  We exploit this fact to present a
  methodology for structural identifiability analysis that detects such
  parameters and characterizes the symmetries in which they are involved.
  We demonstrate the use of our methodology by applying it to four case
  studies.
\end{abstract}

\section{Introduction}

In this paper, we study two closely related properties of nonlinear dynamic
models, namely structural identifiability and observability.  A state
variable is \textit{observable} if its initial value can be determined from
knowledge of the model's subsequent input and output trajectories, and a
parameter in a model is structurally \textit{identifiable} if it can be
determined in the same way. These properties are \textit{structural}, in
the sense that they only depend on the model equations; they are not
affected by the quality of the experimental data.  Since a parameter can be
considered as a constant state variable, structural identifiability will be
treated as a particular case of observability
\cite{chatzis2015observability,villaverde2019observability} and we will
speak of \emph{SIO analysis} of a system.

Nonlinear symmetries \cite{bk:sym,olv:lgde} are useful for analyzing
control systems
\cite{mv:lie,merkt2015higher,skz:lie,sedoglavic2002probabilistic,%
  villaverde2022symmetries,villaverde2019observability,%
  villaverde2021testing,yates2009structural,kz:diss,zs:catao}.  In
particular, the existence of continuous symmetries leaving the output
invariant represents an obstruction to local observability.  Discrete
symmetries have been studied much less than their continuous counterpart.
The most popular approach is due to Hydon \cite{peh:discrete} and an
alternative approach was proposed by Gaeta and Rodr{\'\i}guez
\cite{gr:discrete}.  Both approaches are not completely algorithmic.
Furthermore, Hydon's method can only be used if continuous symmetries
exist. In the particularly interesting situation that a control system is
structurally locally identifiable, but not globally, this will not be the
case.

Recently, two of the present authors analysed for a larger number of
parametric models the difference between \emph{local} and \emph{global}
structural identifiability \cite{rey2023structural}.  If a parameter is
structurally locally identifiable, but not globally (which we denote with
the acronym SLING), then for a given output, finitely many different values
are possible for this parameter.  One goal of the current paper is to argue
that this is related to the existence of \emph{discrete} symmetries.  Based
on this fact, we will present a methodology for distinguishing between
SLING and globally identifiable parameters by searching for discrete
symmetries.

To this end, we will follow an approach proposed by Reid \emph{et al.}\
\cite{rww:discrete} and work with the finite determining system instead of
the infinitesimal one. This system is generally highly nonlinear and thus
hard to analyse.  Reid \emph{et al.}\ relied on an algorithm
\texttt{rifsimp} \cite{rwb:invol} which uses many heuristics and hence may
fail.  We will use rigorous techniques from differential algebra (see
\cite{wys:rkdp} for a concise introduction), in particular the Thomas
decomposition originally developed by Thomas \cite{th:ds,th:sr} and later
revived by Gerdt \cite{vpg:decomp}.  It is implemented in \textsc{Maple} as
a package called \textsc{TDDS} \cite{glhr:tdds} (for an earlier version see
\cite{bglr:thomas}).

It should be noted that, for just deciding SIO, it is not necessary to
\emph{compute} symmetries; it suffices to \emph{detect} their existence.
This can be done algorithmically with methods from the formal theory of
differential equations (see e.\,g.\ \cite{wms:invol} and references therein
for an extensive introduction); some applications of such methods in the
context of symmetry theory have e.\,g.\ been discussed in
\cite{wms:aci,wms:sym,wms:mcm}.

The remainder of this paper is organized as follows: in Section
\ref{sec:discr-cont} we provide background for the study of symmetries and
their relation with structural identifiability and observability; in
Section \ref{sec:fads} we describe our procedure for an SIO analysis, which
is then illustrated with four case studies in Section \ref{sec:ex}; lastly,
we provide some concluding remarks in Section \ref{sec:conclusions}.

\section{Discrete and Continuous Symmetries for Structural Identifiability
  and Observability}\label{sec:discr-cont}

We will discuss finite-dimensional nonlinear control systems of the general form
\begin{subequations}\label{eq:pcsys}
  \begin{align}
    \dot{\mathbf{x}} &=
                       \mathbf{f}(t,\mathbf{x},\mathbf{u},\boldsymbol{\theta})\,,
                       \label{eq:pstate}\\
    \mathbf{y} &=
                 \mathbf{h}(t,\mathbf{x},\mathbf{u},\boldsymbol{\theta})\,.
                 \label{eq:poutput}
  \end{align}
  Here, as usual, $\mathbf{x}\in\RR^{n}$ denotes the \emph{state},
  $\mathbf{u}\in\RR^{m}$ the \emph{input}, $\mathbf{y}\in\RR^{\ell}$ the
  \emph{output} and all these dependent variables are assumed to be smooth
  functions of a single independent variable, the time $t$.  In addition,
  the system depends on parameters $\boldsymbol{\theta}\in\RR^{k}$.  We
  will assume that the right hand sides $\mathbf{f}$ and $\mathbf{h}$ are
  rational functions of their arguments.  This covers most systems
  appearing in typical applications.  For treating structural
  identifiability and observability in a uniform manner, we treat the
  parameters as additional state variables with a trivial dynamics, i.\,e.\
  if we augment \eqref{eq:pcsys} by
  \begin{equation}\label{eq:pdyn}
    \dot{\boldsymbol{\theta}}=0\;.
  \end{equation}
\end{subequations}

A state variable $x_{i}$ is \emph{globally observable} and a parameter
$\theta_{i}$ is \emph{structurally globally identifiable}, respectively, if
for all admissible inputs $\mathbf{u}$, almost all state vectors
$\mathbf{x}$, $\mathbf{x}^{*}$ and almost all parameter vectors
$\boldsymbol{\theta}$, $\boldsymbol{\theta}^{*}$ the equality
$\mathbf{y}(t,\mathbf{x},\mathbf{u},\boldsymbol{\theta})=
\mathbf{y}(t,\mathbf{x}^{*},\mathbf{u},\boldsymbol{\theta}^{*})$ implies
that $x_{i}=x_{i}^{*}$ and $\theta_{i}=\theta_{i}^{*}$, respectively.  A
model is \emph{globally SIO}, if all state variables are globally
observable and all parameters are structurally globally identifiable.  In
the \emph{local} case, the above inequality admits for fixed $x_{i}$ and
$\theta_{i}$, respectively, only finitely many solutions $x_{i}^{*}$ and
$\theta_{i}^{*}$, respectively.

The existence of symmetries among variables may prevent a model from being
fully identifiable and/or observable.  Hence, a goal of this paper is to
devise a method for finding such symmetry transformations.  For an SIO
analysis, one considers nonlinear \emph{point transformations} of the
restricted form\footnote{One could admit that $\mathbf{X}$ and
  $\boldsymbol{\Theta}$ also depends on~$\mathbf{u}$.  But one will find
  among the determining equations
  $\mathbf{X}_{\mathbf{u}}=\boldsymbol{\Theta}_{\mathbf{u}}=0$.  Hence our
  restricted ansatz suffices.}
\begin{subequations}\label{eq:prptrafo}
  \begin{equation}\label{eq:ptrafoxth}
    \tilde{t}=t\,,\quad
    \tilde{\mathbf{x}}=\mathbf{X}(t,\mathbf{x},\boldsymbol{\theta})\,,\quad
    \tilde{\boldsymbol{\theta}}=
    \boldsymbol{\Theta}(t,\mathbf{x},\boldsymbol{\theta})\,,\quad
    \tilde{\mathbf{u}}=\mathbf{u}
  \end{equation}
  Since the input $\mathbf{u}$ is considered as known, we should indeed
  admit only transformations that do not change it.  And as SIO is
  concerned with reconstructing the state $\mathbf{x}$ and the parameters
  $\boldsymbol{\theta}$ from the output $\mathbf{y}$, we cannot allow for
  transformation of the time~$t$.  Finally, the transformation of the
  output $\mathbf{y}$ is determined by entering \eqref{eq:ptrafoxth} into
  \eqref{eq:poutput}.  A point transformation must be invertible requiring
  that the Jacobian
  $\partial(\mathbf{X},\boldsymbol{\Theta})/\partial(\mathbf{x},\boldsymbol{\theta})$
  must possess full rank almost everywhere.
  
  The transformation of derivatives of arbitrary order is completely
  determined by the chain rule (see e.\,g.\ \cite{bk:sym}); for first-order
  derivatives one obtains the \emph{prolongation}
  \begin{align}\label{eq:ptrafodotxth}
    \dot{\tilde{\mathbf{x}}}&=
                              \mathbf{X}^{(1)}(t,\mathbf{x},\boldsymbol{\theta},
                              \dot{\mathbf{x}},\dot{\boldsymbol{\theta}})=
                              \mathbf{X}_{t}+\mathbf{X}_{\mathbf{x}}\dot{\mathbf{x}}+
                              \mathbf{X}_{\boldsymbol{\theta}}\dot{\boldsymbol{\theta}}\,,\\
    \dot{\tilde{\boldsymbol{\theta}}}&=
                                       \boldsymbol{\Theta}^{(1)}(t,\mathbf{x},\boldsymbol{\theta},
                                       \dot{\mathbf{x}},\dot{\boldsymbol{\theta}})=
                                       \boldsymbol{\Theta}_{t}+
                                       \boldsymbol{\Theta}_{\mathbf{x}}\dot{\mathbf{x}}+
                                       \boldsymbol{\Theta}_{\boldsymbol{\theta}}\dot{\boldsymbol{\theta}}\,.
  \end{align}
\end{subequations}

A control system \eqref{eq:pcsys} is \emph{not} globally SIO, if it
possesses symmetries that leave the output invariant.  Hence, we require
that both the state equation \eqref{eq:pstate} and the output equation
\eqref{eq:poutput} remain invariant under the point transformation
\eqref{eq:prptrafo}.  Entering \eqref{eq:prptrafo} into \eqref{eq:pcsys}
yields the conditions
\begin{subequations}\label{eq:ptrsys}
  \begin{gather}
    \mathbf{X}^{(1)}(t,\mathbf{x},\boldsymbol{\theta},
    \dot{\mathbf{x}},\dot{\boldsymbol{\theta}})=
    \mathbf{f}\bigl(t,\mathbf{X}(t,\mathbf{x},\boldsymbol{\theta}),
    \mathbf{u},\boldsymbol{\Theta}(t,\mathbf{x},\boldsymbol{\theta})\bigr)\,,\\
    \boldsymbol{\Theta}^{(1)}(t,\mathbf{x},\boldsymbol{\theta},
    \dot{\mathbf{x}},\dot{\boldsymbol{\theta}})=0\,,\\
    \mathbf{h}(t,\mathbf{x},\mathbf{u},\boldsymbol{\theta})=
    \mathbf{h}\bigl(t,\mathbf{X}(t,\mathbf{x},\boldsymbol{\theta}),
    \mathbf{u},\boldsymbol{\Theta}(t,\mathbf{x},\boldsymbol{\theta})\bigr)\,.
    \label{eq:ptrsyso}
  \end{gather}
\end{subequations}
Note that \eqref{eq:ptrsyso} is a purely algebraic equation, as no
derivatives appear in it.

In these equations, we eliminate $\dot{\mathbf{x}}$ and
$\dot{\boldsymbol{\theta}}$, respectively, by substituting them by
$\mathbf{f}(t,\mathbf{x},\mathbf{u},\boldsymbol{\theta})$ and~$0$,
respectively, i.\,e.\ by using the state equations.  Then we clear
denominators (this is not necessary, if \eqref{eq:pcsys} is a polynomial
system) to obtain expressions that are multivariate polynomials in the
inputs $\mathbf{u},$ with coefficients $\Delta_{j}$ which are polynomials
in the variables $t$, $\mathbf{x}$, $\boldsymbol{\theta}$ and the functions
$\mathbf{X}$, $\boldsymbol{\Theta}$ plus their derivatives
$\mathbf{X}_{t}$, $\mathbf{X}_{\mathbf{x}}$, $\boldsymbol{\Theta}_{t}$,
$\boldsymbol{\Theta}_{\mathbf{x}}$ (since we have for the parameters the
trivial dynamics $\dot{\boldsymbol{\theta}}=0$, the derivatives
$\mathbf{X}_{\boldsymbol{\theta}}$ and
$\boldsymbol{\Theta}_{\boldsymbol{\theta}}$ will not show up).  The
\emph{finite determining system} $\findet$ for point symmetries of
\eqref{eq:pcsys} arises by requiring that all these coefficients
$\Delta_{j}$ vanish, as the invariance must hold for arbitrary
inputs~$\mathbf{u}$.  It represents a complicated system of first-order
polynomially nonlinear partial differential equations,
\begin{equation}\label{eq:rfindetsys}
  \Delta_{j}(t,\mathbf{x},\boldsymbol{\theta},\mathbf{X},\boldsymbol{\Theta},
  \mathbf{X}_{t},\boldsymbol{\Theta}_{t}
  ,\mathbf{X}_{\mathbf{x}},\boldsymbol{\Theta}_{\mathbf{x}})=0\,,
  j=1,\dots,J\,,
\end{equation}
for the functions $\mathbf{X}$ and $\boldsymbol{\Theta}$ which is generally
not solved for any derivatives, i.\,e.\ fully implicit.  The number $J$ of
equations in \eqref{eq:rfindetsys} depends not only on the dimensions $n$,
$m$, $\ell$, $k$ of the vectors $\mathbf{x}$, $\mathbf{u}$, $\mathbf{y}$,
$\boldsymbol{\theta}$, but also on the degrees with which the inputs
$\mathbf{u}$ appear in \eqref{eq:pcsys}.  If we assume that
\eqref{eq:pcsys} is affine in the inputs (as is often the case in
applications), then $J=(n+\ell+k)(m+1)$.  The solution set of $\findet$ has
the structure of a local group $\mathcal{G}$, i.\,e.,\ where defined, the
composition of two solutions defines again a solution, the inverse of a
solution is also a solution, and the identity map is a solution.

Most people using symmetry methods do not set up $\findet$.  Instead, they
use the infinitesimal approach developed by Lie already in the 19th
century.  Mathematically, it corresponds to determining the Lie
algebra~$\mathfrak{g}$ of the Lie group~$\mathcal{G}$ of continuous point
symmetries.  In this approach, one assumes that one has at least a
one-parameter group of transformations, i.\,e.\ one must augment
\eqref{eq:ptrafoxth} by a parameter~$\epsilon$:
\begin{equation}\label{eq:oneparaxth}
  \tilde{t}=t\,,\ \
  \tilde{\mathbf{x}}=\mathbf{X}(t,\mathbf{x},\boldsymbol{\theta},\epsilon)\,,\ \
  \tilde{\boldsymbol{\theta}}=
    \boldsymbol{\Theta}(t,\mathbf{x},\boldsymbol{\theta},\epsilon)\,,\ \ \tilde{\mathbf{u}}=u\,.
\end{equation}
Here, we always assume that we obtain the identity transformation for
$\epsilon=0,$ and that the concatenation of two transformations with the
parameter values $\epsilon_{1}$ and $\epsilon_{2}$, respectively, yields
the transformation for $\epsilon_{1}+\epsilon_{2}$.

For small $\epsilon$, one can linearise \eqref{eq:oneparaxth} around
$\epsilon=0$ and obtains the corresponding \emph{infinitesimal generator},
the vector field
\begin{equation}\label{eq:infgen}
  V=
  \frac{\partial\mathbf{X}}{\partial\epsilon}\bigg|_{\epsilon=0}
  \!\!\!\cdot\partial_{\mathbf{x}}  +
 \frac{\partial\boldsymbol{\Theta}}{\partial\epsilon}\bigg|_{\epsilon=0}
  \!\!\!\cdot\partial_{\boldsymbol{\theta}} \,.
\end{equation}
One can obtain the \emph{prolonged infinitesimal generator} $V^{(1)}$ by
linearising \eqref{eq:ptrafodotxth}.  But given an arbitrary vector field
\begin{equation}\label{eq:ransatz}
  V
  =\boldsymbol{\xi}(t,\mathbf{x},\boldsymbol{\theta})
  \cdot\partial_{\mathbf{x}}+
  \boldsymbol{\zeta}(t,\mathbf{x},\boldsymbol{\theta})
  \cdot\partial_{\boldsymbol{\theta}}\,,
\end{equation}
$V^{(1)}$ can also be determined directly.  The ansatz
\begin{equation}
  V^{(1)}=V+
  \boldsymbol{\xi}^{(1)}(t,\mathbf{x},\boldsymbol{\theta},
  \dot{\mathbf{x}},\dot{\boldsymbol{\theta}})
  \cdot\partial_{\dot{\mathbf{x}}}+
  \boldsymbol{\zeta}^{(1)}(t,\mathbf{x},\boldsymbol{\theta},
  \dot{\mathbf{x}},\dot{\boldsymbol{\theta}})
  \cdot\partial_{\dot{\boldsymbol{\theta}}}
\end{equation}
yields after a short calculation the following explicit representation for
the coefficients
\begin{subequations}
  \begin{align}
    \boldsymbol{\xi}^{(1)} &=
                             \boldsymbol{\xi}_{t} +
                             \boldsymbol{\xi}_{\mathbf{x}}\dot{\mathbf{x}} +
                             \boldsymbol{\xi}_{\boldsymbol{\theta}}\dot{\boldsymbol{\theta}}\,,\\
    \boldsymbol{\zeta}^{(1)} &=
                               \boldsymbol{\zeta}_{t} +
                               \boldsymbol{\zeta}_{\mathbf{x}}\dot{\mathbf{x}} +
    \boldsymbol{\zeta}_{\boldsymbol{\theta}}\dot{\boldsymbol{\theta}}\,.
  \end{align}
\end{subequations}
These formulae can be extended to higher derivatives \cite{bk:sym}, but the
expressions are getting rapidly rather complicated.

The ansatz \eqref{eq:ransatz} defines an \emph{infinitesimal symmetry} of
the control system \eqref{eq:pcsys}, if and only if it satisfies
\begin{subequations}\label{eq:infsymmpout}
  \begin{gather}
    V^{(1)}\bigl(\dot{\mathbf{x}}-\mathbf{f}(t,\mathbf{x},\mathbf{u},\boldsymbol{\theta})\bigr)
    _{\raisebox{0.5ex}{\big|}\dot{\mathbf{x}}=\mathbf{f}(t,\mathbf{x},\mathbf{u},\boldsymbol{\theta}),\
      \dot{\boldsymbol{\theta}}=0}=0\,,\\
    V^{(1)}\bigl(\dot{\boldsymbol{\theta}}\bigr)
    _{\raisebox{0.5ex}{\big|}\dot{\mathbf{x}}=\mathbf{f}(t,\mathbf{x},\mathbf{u},\boldsymbol{\theta}),\
      \dot{\boldsymbol{\theta}}=0}=0\,,\\
    V\mathbf{h}(t,\mathbf{x},\mathbf{u},\boldsymbol{\theta})=0\,.
  \end{gather}
\end{subequations}
The left-hand side of the first equation means that one applies first the
prolonged generator $V^{(1)}$ to the expression
$\dot{\mathbf{x}}-\mathbf{f}(t,\mathbf{x},\mathbf{u})$ and then replaces in
the result every occurrence of $\dot{\mathbf{x}}$ by
$\mathbf{f}(t,\mathbf{x},\mathbf{u})$.  For a rational $\mathbf{f}$, one
has to clear denominators in the result.  The second equation is easy to
evaluate: it yields the partial differential equation
\begin{equation}
  \boldsymbol{\zeta}^{(1)}\bigl(t,\mathbf{x},
  \mathbf{f}(t,\mathbf{x},\mathbf{u},\boldsymbol{\theta}),0\bigr)=
  \boldsymbol{\zeta}_{t} +
  \boldsymbol{\zeta}_{\mathbf{x}}\cdot\mathbf{f}(t,\mathbf{x},\mathbf{u},\boldsymbol{\theta})=0\,.
\end{equation}
Now one can again consider each equation in \eqref{eq:infsymmpout} as a
multivariate polynomial in the inputs~$\mathbf{u}$ with coefficients
$\delta_{j}$ which are polynomials in the variables $t$, $\mathbf{x}$,
$\boldsymbol{\theta}$ and the functions $\boldsymbol{\xi}$,
$\boldsymbol{\zeta}$ plus their partial derivatives $\boldsymbol{\xi}_{t}$,
$\boldsymbol{\xi}_{\mathbf{x}}$, $\boldsymbol{\xi}_{\boldsymbol{\theta}}$,
$\boldsymbol{\zeta}_{t}$, $\boldsymbol{\zeta}_{\mathbf{x}}$,
$\boldsymbol{\zeta}_{\boldsymbol{\theta}}$ and the \emph{infinitesimal
  determining system} $\infdet$ is obtained by requiring that all
coefficients of these polynomials vanish.  In contrast to the finite
determining system, we obtain now \emph{linear} partial differential
equations
\begin{equation}\label{eq:infdetsys}
  \delta_{j}(t,\mathbf{x}, \boldsymbol{\theta},\boldsymbol{\xi},\boldsymbol{\zeta},
  \boldsymbol{\xi}_{t},\boldsymbol{\zeta}_{t},\boldsymbol{\xi}_{\mathbf{x}},
  \boldsymbol{\zeta}_{\mathbf{x}},\boldsymbol{\xi}_{ \boldsymbol{\theta}},
  \boldsymbol{\zeta}_{ \boldsymbol{\theta}})=0\,,
  j=1,\dots,J\,,
\end{equation}
for the coefficients $\boldsymbol{\xi}$, $\boldsymbol{\zeta}$ of the vector
field~$V$.

The solution space of $\infdet$ defines a Lie algebra~$\mathfrak{g}$,
i.\,e.\ if we consider solutions as vector fields, then the Lie bracket of
two solutions is again a solution.  Furthermore, one can obtain from
$\mathfrak{g}$ a connected Lie group $\mathcal{G}_{0}$ by
\emph{exponentiation}: if
$\hat{\boldsymbol{\xi}}(t,\mathbf{x}, \boldsymbol{\theta})$ and
$\hat{\boldsymbol{\zeta}}(t,\mathbf{x}, \boldsymbol{\theta})$ form a
solution of $\infdet$, then integrating the ordinary differential equations
$\frac{d\mathbf{X}}{d\epsilon}=\hat{\boldsymbol{\xi}}(T,\mathbf{X},
\boldsymbol{\Theta})$,
$\frac{d\boldsymbol{\Theta}}{d\epsilon}=\hat{\boldsymbol{\zeta}}(T,\mathbf{X},
\boldsymbol{\Theta})$ yields a one-parameter group of transformations of
the form \eqref{eq:oneparaxth}.  For any choice of the parameter
$\epsilon$, the obtained point transformation will solve the finite
determining system~$\findet$.  But we cannot expect that \emph{all}
solutions of $\findet$ arise in this way, i.\,e.\ that they can be embedded
in a one-parameter family of transformations.  In other words: generally
$\mathcal{G}$ is larger than $\mathcal{G}_{0}$.

From a mathematical point of view, the deeper reason is that one assumes in
the infinitesimal approach that $\mathcal{G}$ contains a whole, non-trivial
(local) \emph{Lie group} of symmetries, i.\,e.\ \emph{continuous
  symmetries}.  Even if this is the case, this group may consist of several
disjoint components and infinitesimal methods see only the connected
component~$\mathcal{G}_{0}$ of the identity.  If a dynamical system has
only \emph{discrete symmetries}, then $\mathcal{G}$ is not a Lie group, but
typically a finite group, and not visible for infinitesimal techniques.

If a control system \eqref{eq:pcsys} admits continuous symmetries
preserving the output, then it cannot be locally SIO, as infinitely many
solutions yield the same output.  These solutions form a continuous family
and thus cannot even locally be distinguished. The absence of such
symmetries guarantees local SIO \cite[Cor.~3.1]{yates2009structural}.  In
the case of a finite group of discrete symmetries, only finitely many
solutions lead to the same output.  The system is then locally SIO, as in a sufficiently small neighbourhood only one appropriate solution exists, but not globally.  

\section{Formal Analysis of Determining Systems}
\label{sec:fads}

Even the linear determining system $\infdet$ for \emph{infinitesimal}
symmetries is often difficult to solve explicitly for systems of ordinary
differential equations (the situation is fairly different for partial
differential equations where such determining systems are routinely solved
by computer algebra packages relying on a combination of heuristics and
systematic theory).  The \emph{finite} determining system $\findet$ will be
explicitly solvable only in exceptional cases.  However, a \emph{formal
  analysis} of the determining systems is possible much more
often. ``Formal'' has here two meanings: it implies that firstly, one only
works with the differential equations themselves, trying, for example, to
bring them into a suitable normal form, and that secondly, one considers
only formal power series solutions (see \cite{wms:invol} for a more
extensive discussion).

In our context, one of the key quantities is the size of the solution
space: if more than one solution exists, then non-trivial symmetries exist,
and the given system is not SIO.  There are methods for formally counting
solutions (or more precisely freely choosable Taylor coefficients of formal
solutions) --- see e.\,g.\ \cite{lh:phd,lh:dcp,wms:str} --- and for the
infinitesimal determining system $\infdet$ such a counting was already
performed by Schwarz \cite{fs:size}.

In this paper, we concentrate only on distinguishing between infinitely
many and only finitely many solutions corresponding to the distinction
between control systems which are not even locally SIO and those that are
locally, but not globally SIO.  For this purpose, it suffices to bring the
determining system into a kind of normal form from which the size of the
solution space can be read off. For the linear infinitesimal determining
system $\infdet$, \emph{Janet bases} provide such a normal form
\cite{wms:invol,sz:atls}; in the nonlinear case, one needs the \emph{Thomas
  decomposition} \cite{dr:habil} for analysing $\findet$.

As $\infdet$ always has the zero solution and $\findet$ the identity map as
trivial solutions, determining systems can never be inconsistent.  If
$\infdet$ has more solutions than only the zero solution, then it has
automatically infinitely many solutions implying the studied control system
\eqref{eq:pcsys} is not even locally SIO.  This is easy to detect from a
Janet basis, as the basis becomes trivial, if only the zero solution
exists.

The analysis of nonlinear systems is more involved. The first difference is
that one cannot simply transform $\findet$ into a normal form, but must
perform case distinctions so that the decomposition actually consists of a
finite number of so-called \emph{simple systems} \cite{dr:habil} which have
to be studied separately.  We have infinitely many solutions, if at least
one of these simple systems contains an equation which is still a
differential equation, i.\,e.\ in which still a derivative appears. Indeed,
in this case, one can choose freely at least one Taylor coefficient of the
formal solution.

The simplest case for a finite solution space arises if each simple system
admits exactly one solution (in this case, the simple system contains for
each unknown function one linear equation determining it uniquely).  Then
the total number of solutions is just the number of simple systems.
Generally, a simple system may also contain algebraic equations of higher
degree.  If we define the degree of a simple system as the product of the
degrees of all the equations contained in it, then the defining properties
of a simple system ensure that the number of solutions is exactly this
degree.  The total number of solutions of the determining system is then
the sum of the degrees of the arising simple systems, as the theory of the
Thomas decomposition asserts that the solution spaces of the simple systems
are disjoint (this is an important difference to an alternative
decomposition method due to Boulier \emph{et al.}\
\cite{blop:radical,blop:radical2} where solution spaces can intersect).

\begin{figure}[t]
  \centering
\begin{tikzpicture}[-latex,thick,font=\small]
  
  \node[draw,
    trapezium, 
    trapezium left angle = 65,
    trapezium right angle = 115,
    trapezium stretches,
    top color=white, bottom color=green!30,
    minimum width = 7em, minimum height = 5ex, inner sep=1ex
    ] (input) {\textbf{control system}};

  \node[draw,
    shape=rectangle, rounded corners,
    below left=of input,
    top color=white, bottom color=blue!20,
    minimum width = 7em, minimum height = 5ex, inner sep=1ex
    ] (block11) {set up $\Delta_{\mathrm{fin}}$};

  \node[draw,
    shape=rectangle, rounded corners,
    below right=of input,
    top color=white, bottom color=blue!20,
    minimum width = 7em, minimum height = 5ex, inner sep=1ex
    ] (block12) {set up $\Delta_{\mathrm{inf}}$};

  \node[draw,
    shape=rectangle, rounded corners,
    below=of block11, align=center,
    top color=white, bottom color=blue!20,
    minimum width = 7em, minimum height = 5ex, inner sep=1ex
    ] (block21) {compute$\strut$ Thomas\\ decomposition$\strut$};

  \node[draw,
    shape=rectangle, rounded corners,
    below=of block12, align=center,
    top color=white, bottom color=blue!20,
    minimum width = 7em, minimum height = 5ex, inner sep=1ex
    ] (block22) {compute$\strut$\\ Janet basis$\strut$};

  \node[draw,
    shape=rectangle, rounded corners,
    below=of block21, align=center,
    top color=white, bottom color=blue!20,
    minimum width = 7em, minimum height = 5ex, inner sep=1ex
    ] (block31) {determine\\ size $S$ of ${\cal G}$};

  \node[draw,
    shape=rectangle, rounded corners,
    below=of block22, align=center,
    top color=white, bottom color=blue!20,
    minimum width = 7em, minimum height = 5ex, inner sep=1ex
    ] (block32) {determine size\\ $S_{0}$ of ${\cal G}_{0}$};

  \node[draw,
    diamond,
    below=of $(block31)!0.5!(block32)$,
    top color=white, bottom color=red!30,
    minimum width = 5em, inner sep=0ex
    ] (check) {$S>S_{0}$};
    
  \node[draw,
    trapezium, 
    trapezium left angle = 65,
    trapezium right angle = 115,
    trapezium stretches,
    below left=of check, xshift=-2em, align=center,
    top color=white, bottom color=green!30,
    minimum width = 7em, minimum height = 5ex, inner sep=1ex
    ] (outputY) {\textbf{discrete}\\ \textbf{symmetries}};

  \node[draw,
    trapezium, 
    trapezium left angle = 65,
    trapezium right angle = 115,
    trapezium stretches,
    below right=of check, xshift=2em, align=center,
    top color=white, bottom color=green!30,
    minimum width = 7em, minimum height = 5ex, inner sep=1ex
    ] (outputN) {\textbf{no discrete}\\ \textbf{symmetries}};

    \draw (block11) edge (block21);
    \draw (block21) edge (block31);
    \draw (block12) edge (block22);
    \draw (block22) edge (block32);
    \draw (input) |- (block11);
    \draw (input) |- (block12);
    \draw (block31) -| (check);
    \draw (block32) -| (check);
    \draw (check) -| (outputY) node[pos=0.75,fill=white,inner sep=3pt]{yes};
    \draw (check) -| (outputN) node[pos=0.75,fill=white,inner sep=3pt]{no};
\end{tikzpicture}
\caption{Basic workflow of our approach to detecting the existence of
  discrete symmetries and, as a result, establishing that a model is
  locally but not globally identifiable.}
\label{fig:flow}
\end{figure}

The basic idea of our approach to detect discrete symmetries consists of
comparing the results of a formal analysis of $\findet$ with those of a
formal analysis of $\infdet$.  A sketch is given in Figure~\ref{fig:flow},
where the left column depicts the analysis of $\findet$ and the right
column the analysis of $\infdet$.  For lack of space, we cannot discuss
here how one can actually measure the size of a solution space with the
so-called \emph{differential counting polynomial} \cite{lh:phd,lh:dcp}.  We
note that while the solution space of $\infdet$ is the Lie algebra
$\mathfrak{g}$ and not the Lie group $\mathcal{G}_{0}$, the sizes of
$\mathfrak{g}$ and $\mathcal{G}_{0}$ are identical, as $\mathcal{G}_{0}$
can -- in principle -- be computed from $\mathfrak{g}$ by exponentiation,
i.\,e.\ the solution of a system of ordinary differential equations, which
defines a bijection between the two sets.  For the purposes of this paper,
we are mainly interested in the case that $\mathcal{G}_{0}$ is the trivial
group consisting only of the identity.  If now $\mathcal{G}$ contains more
elements, then we have only discrete symmetries and thus are dealing with a
control system which is locally SIO, but not globally.

\section{Examples}
\label{sec:ex}

We present four different examples. In the first one, a simple linear
compartmental model from physiology, only discrete symmetries appear (hence
the system is locally SIO).  As they consist of a simple permutation, they
could have been guessed by a direct inspection of the system.  The main
advantage of our systematic approach is that it guarantees that there are
indeed no further symmetries (continuous or discrete). In what follows, if
a simple system does not include an equation for a transformed variable, it
is understood that this variable remains unchanged.

\begin{example}\label{ex:bilirubin}
  The 4-compartmental mammillary model from \cite[Ex.~7]{mkd:combos}
  consists of the linear control system:
  \begin{equation}\label{eq:bilirubin2}
    \begin{cases}
      \dot{x}_{1} = -(k_{21}+k_{31}+k_{41}+k_{01})x_{1} + {}\\
      \qquad\ k_{12}x_{2} + k_{13}x_{3} + k_{14}x_{4} + u\,,\\
      \dot{x}_{2} = k_{21}x_{1} - k_{12}x_{2}\,,\\
      \dot{x}_{3} = k_{31}x_{1}- k_{13}x_{3}\,,\\
      \dot{x}_{4} = k_{41}x_{1}- k_{14}x_{4}\,,\\
      y=x_{1}\,.
    \end{cases}
  \end{equation}
  By direct inspection, one easily identifies an invariance under an action
  of the finite group $S_{3}$ (the symmetric group for sets with three
  elements) containing six elements: if we form the three triples
  $(k_{12},k_{13},k_{14})$, $(k_{21},k_{31},k_{41})$, $(x_{2},x_{3},x_{4})$
  and perform on them simultaneously the same permutation of the indices
  $2,3,4$, then the system \eqref{eq:bilirubin2} remains unchanged. Setting
  up $\findet$ and applying the Thomas decomposition reveals that this is
  indeed the complete solution space, i.\,e.\ that
  $\mathcal{G}=S_{3}$. Each of the six arising simple systems admits
  exactly one solution representing one of the transformations described
  above. One of these simple systems is shown here as a representative
  example:
  \begin{gather*}
    \tilde{x}_2 = x_3,\ \tilde{x}_3 = x_4,\ \tilde{x}_4 = x_2,\ \tilde{k}_{12} = k_{13},\ \tilde{k}_{13} =
    k_{14},\ \\ \tilde{k}_{14} = k_{12},\ \tilde{k}_{21} = k_{31},\ \tilde{k}_{31} = k_{41},\ \tilde{k}_{41} = k_{21}.
  \end{gather*}
  Note that because of its linearity, one can directly read off the
  permutation it defines.  Hence we can conclude that in this example all
  states are locally observable and all parameters are locally
  identifiable.  But only $x_1$ and $k_{01}$ are globally observable and
  identifiable, respectively, as they are not affected by the permutations.
\end{example}

Our second example is a well studied model from biology with a more complex
finite solution space, demonstrating the rise of multiple solutions from a single
simple system.

\begin{example}
  We consider the Goodwin oscillator \cite{goodwin_OscillatoryBehaviorEnzymatic_1965}
  \begin{equation}
    \begin{cases}
      \dot{X} = \frac{1}{1+Z^{m}}-X\,,\\
      \dot{Y} = X - k_2 Y\,,\\
      \dot{Z} = k_1 Y - k_3 Z\,,\\
      y = X
    \end{cases}
  \end{equation}
  in a non-dimensional form. For a Hill coefficient $m=4$, the finite
  determining system decomposes into six simple systems of purely algebraic
  equations. Two of them involve quadratic equations in $k_1$ and the
  remaining systems are linear, thus the discrete symmetries form a group
  of order~$8$. It is not difficult to identify it as the commutative group
  $\mathcal{G}=C_4\times C_2$ (the product of two cyclic groups). The
  generator for the subgroup $C_2$ is given as the single solution of the
  simple system
  \begin{equation}\label{eq:goodC2}
    \tilde{Y}=Y+\frac{Z}{k_1}(k_2-k_3),\quad 
    \tilde{k}_2 = k_3,\quad \tilde{k}_3 = k_2.
  \end{equation}
  It is easy to see that this transformation does not affect the first
  equation of the system nor the output and is a symmetry of the other two
  equations. The existence of this subgroup implies that the state $Y$ is
  only locally observable and that the parameters $k_2$ and $k_3$ are
  SLING.  As a generator of the cyclic subgroup $C_4$, one can take any of
  the two solutions of the simple system
  \begin{equation}
    \tilde{Z}=-\frac{k_1}{\tilde{k}_1}Z,\quad \tilde{k}_1^2 = -k_1^2.
  \end{equation}
  These transformations leave both the Hill function in the first equation
  and the complete third equation invariant and do not affect the remainder
  of the system. However, in contrast to the transformation
  \eqref{eq:goodC2}, they are biologically not relevant, as they lead to
  negative or even complex values for $k_1$ and $Z$.  Thus, from a
  biological point of view, $Z$ is globally observable and $k_1$ globally
  identifiable.  One can verify by direct computation that for any Hill
  coefficient $m\in\mathbbm{N}$ one finds the commutative group
  $C_m\times C_2$ as a discrete symmetry group and we conjecture that, as
  for $m=4$, there are no other symmetries.
\end{example}

Our next example is a classical model from the control theory literature.
Like the first example, it possesses an easy to find permutation symmetry,
but in addition also a continuous symmetry group.

\begin{example}\label{ex:llw1987}
  The following model stems from \cite{llw:1987}:
  \begin{equation}\label{eq:llw}
    \begin{cases}
      \dot{x}_{1} = -\theta_{1}x_{1}+\theta_{2}u\,,\\
      \dot{x}_{2} = -\theta_{3}x_{2}+\theta_{4}u\,,\\
      \dot{x}_{3} = -(\theta_{1}+\theta_{3})x_{3}+
      (\theta_{4}x_{1}+\theta_{2}x_{2})u\,,\\
      y = x_{3}\;.
    \end{cases}
  \end{equation}
  Direct inspection yields again a discrete $S_{2}$ symmetry: if we form
  the three pairs $(\theta_{1},\theta_{3})$, $(\theta_{2},\theta_{4})$,
  $(x_{1},x_{2})$ and simultaneously swap their elements, then
  \eqref{eq:llw} remains unchanged.  However, this time the symmetry group
  $\mathcal{G}$ is larger.  The Thomas decomposition of $\findet$ consists
  of two simple systems, each comprising a set of algebraic equations and
  an identical set of differential equations for $\tilde{\theta}_4$,
  indicating the presence of infinitely many solutions. Solving for
  $\tilde{\theta}_4$ reveals that the infinite part of the solution space
  can be parametrized by a single function of the form
  \begin{displaymath}
      \varphi:=\varphi\Bigl(\theta_1, \theta_2, \theta_3, \theta_4,
      (-x_1x_2 + x_3)e^{t(\theta_1 + \theta_3)}\Bigr)\,.
  \end{displaymath}
  This result implies that the model \eqref{eq:llw} is invariant under a
  Lie group $\mathcal{G}$ consisting of two connected components.  The
  component $\mathcal{G}_{0}$ containing the identity is given by the
  transformations
  \begin{gather*}
    \tilde{x}_1 = \frac{\theta_4x_1}{\varphi}\,,\
    \tilde{x}_2 = \frac{x_2\varphi}{\theta_4}\,,\\
    \tilde{\theta}_1 = \theta_1\,,\
    \tilde{\theta}_2 = \frac{\theta_2\theta_4}{\varphi}\,,\
    \tilde{\theta}_3 = \theta_3\,,\
    \tilde{\theta}_4 = \varphi
  \end{gather*}
  (the identity is obtained by choosing $\varphi=\theta_4$).  The second
  connected component is given by the transformations
  \begin{gather*}
   \quad \tilde{x}_1 = \frac{\theta_2x_2}{\varphi}\,,\
   \tilde{x}_2 = \frac{x_1\varphi}{\theta_2}\,,\\
    \tilde{\theta}_1 = \theta_3\,,\
    \tilde{\theta}_2 = \frac{\theta_2\theta_4}{\varphi}\,,\
    \tilde{\theta}_3 = \theta_1\,,\
    \tilde{\theta}_4 = \varphi
  \end{gather*}
  (the second permutation mentioned above arises by choosing
  $\varphi=\theta_2$).  Thus, the symmetry analysis reveals that the two
  parameters $\theta_1$ and $\theta_3$ are SLING, whereas all other states
  and parameters are not even locally observable or identifiable.
\end{example}

Our last example is an epidemiological model where the appearing discrete
symmetry group $\mathcal{G}$ operates in a highly non-trivial way and where it appears impossible to guess the form of the symmetries just by a simple inspection of the control system.

\begin{example}\label{ex:seirq}
  In \cite{zha:seirq}, the authors augmented a classical SEIR model for an
  epidemic by a quarantine compartment, which is also the output, leading
  to the input-free system
  \begin{equation}\label{eq:seirq}
    \begin{cases}
      \dot{S} = -\beta S I\,,\\
      \dot{E} = \beta S I - \nu E\,,\\
      \dot{I} = \nu E - \psi I - (1-\psi)\gamma I\,,\\
      \dot{Q} = \psi I -\gamma Q\,,\\
      \dot{R} = (1-\psi)\gamma I + \gamma Q\,,\\
      y = Q\,.
    \end{cases}
  \end{equation}
  Here, not only direct inspection fails, but even the application of our approach reaches its (computational) limits. Hence, we applied certain
  simplifications to $\findet$. The dynamics of the state $R$ in
  \eqref{eq:seirq} is independent of the rest of the system and does not
  affect the output. An analysis of the infinitesimal determining system
  $\infdet$ reveals that only $R$ admits continuous symmetries.
  Considering only the reduced system, obtained by removing the equation
  for $R$, simplifies $\findet$. Furthermore, we adopt our heuristic that
  the transformed parameters are independent of~$\mathbf{x}$ and assume
  that $\gamma$ remains invariant. The resulting determining system then
  decomposes into two simple systems: one corresponding to the identity
  transformation and the other yielding a non-trivial symmetry:
  \begin{align*}
    \tilde{S} &= \frac{\nu \psi S (\gamma - 1)}{(\nu - \gamma)(\psi\gamma - \psi -
    \gamma)},\ \tilde{I} = \frac{I\psi(1-\gamma)}{\nu - \gamma}\,, \\
    \tilde{E} &= \frac{\bigl(I(\psi - 1)\gamma - \psi I + \nu(E + I)\bigr)(\gamma - 1)\psi}{(\nu -
    \gamma)\bigl((\psi - 1)\gamma - \psi\bigr)}\,, \\
    \tilde{\beta} &= \frac{\beta(\gamma-\nu)}{(\gamma - 1)\psi},\ \tilde{\nu} = (1 -
    \psi)\gamma + \psi,\ \tilde{\psi} = \frac{\gamma-\nu}{\gamma - 1}\,.
  \end{align*}
  Due to the made simplifications, we cannot exclude the existence of
  further discrete symmetry transformations. Nevertheless, we can make
  statements about the identifiability and observability. While the state
  $R$ is not even locally observable, the states $S$, $E$, $I$ are locally
  but not globally observable. In addition, the parameters
  $\beta, \nu, \psi$ are SLING.  No statement about $\gamma$ is possible.
  The explicit form of the non-trivial symmetry allows us furthermore to
  conclude that the found discrete symmetry is biologically not relevant:
  biologically meaningful values of $\gamma$ and $\nu$ lie between $0$ and
  $1$ and normally $\gamma$ is larger than $\nu$; applying the above
  transformation then leads to negative values of $\tilde{\psi}$ and
  $\tilde{I}$ which are biologically meaningless.
\end{example}

In the discussion of these examples, $\infdet$ was only mentioned in the
last one, which may appear strange in view of our flowchart in
Figure~\ref{fig:flow}.  The reason is simply that the first three examples
were so ``nice'' that one could explicitly solve $\findet$.  Thus, it was
not necessary to formally count its solutions and to compare with
$\infdet$; we could directly discuss the arising point transformations.
The last example indicates that this is not always the case and there we
needed results from an analysis of $\infdet$ to find useful
simplifications.

\section{Conclusions}\label{sec:conclusions}

In this paper, we have presented an approach to study the structural
identifiability and observability (SIO) of nonlinear dynamic models using
discrete symmetries. It is known that continuous symmetries are sources of
non-identifiability and/or non-observability, i.e., a single input-output
trajectory is compatible with an infinite number of parameter or state
variable values. This fact has already been exploited by several authors
\cite{yates2009structural,merkt2015higher,anguelova2012minimal,meshkat2014identifiable,castro2020testing},
who introduced methods to analyse SIO by searching for continuous
symmetries. However, said methods can only distinguish between
non-identifiability and (at least) local identifiability (and similarly
observability).  In contrast, here we have exploited the fact that the
existence of purely discrete symmetries does not lead to non-identifiable
parameters, but to parameters that are structurally locally identifiable,
but not globally (SLING). Thus, by determining the discrete symmetries in a
model we obtain a more precise SIO analysis, which can distinguish between
SLING and globally identifiable parameters. Here, we have presented a
methodology to perform this analysis. The code that implements the method
and reproduces the results can be accessed at:
\url{https://doi.org/10.5281/zenodo.16410656}.  Key elements in our
procedure are obtaining the finite determining system, analyzing it using
the Thomas decomposition, and finding the number of solutions.  Our method
provides more information than differential algebraic algorithms like the
ones described in \cite{hong2019sian,dong2023differential}, since it can
fully characterize the form of the symmetries in which the parameters and
state variables are involved.  This information is valuable for the
purposes of finding all possible solutions and reparameterizing the model
\cite{meshkat2014identifiable,massonis2023autorepar}.  A future line of
work is to improve the computational efficiency of our method.
Additionally, we envision its extension to systems of partial differential
equations, for which comparatively few techniques for structural
identifiability analysis exist.

\section*{Acknowledgements}

The work of NB, MAM and WMS was performed within the Research Training
Group \emph{Biological Clocks on Multiple Time Scales} (GRK 2749) at Kassel
University funded by the German Research Foundation (DFG). The work of XRB
and AFV was supported by grant PID2023-146275NB-C21 funded by
MICIU/AEI/10.13039/501100011033 and ERDF/E, grant CNS2023-144886 funded by
MICIU/AEI/10.13039/501100011033 and the European Union
Next\-GenerationEU/PRTR, grant RYC-2019-027537-I funded by
MCIN/AEI/10.13039/ 501100011033 and by ``ESF Investing in your future'',
and grant ED431F 2021/ 003 funded by the Xunta de Galicia, Conseller\'ia differential equation
Cultura, Educaci\'on e Universidade. XRB was also partially funded by a PhD
scholarship from the Universidade de Vigo.

\bibliography{IdentSymm}
\bibliographystyle{plain}

\end{document}